\documentclass{article}
\usepackage[utf8]{inputenc}
\usepackage[square,sort,comma,numbers]{natbib}
\usepackage[english]{babel}
\usepackage{graphicx}
\usepackage{hyperref}
\usepackage{amsmath}
\usepackage{amsthm}
\usepackage{amssymb}
\usepackage{dsfont}
\usepackage{authblk}
\newtheorem{theorem}{Theorem}[section]
\newtheorem{conjecture}[theorem]{Conjecture}
\usepackage[flushleft]{threeparttable}
\usepackage{color,soul}
\usepackage{comment}

\usepackage[msc-links]{amsrefs}
\begin{document}

\title{Patterns in Knot  Floer Homology}

\author[1]{Ekaterina S. Ivshina}
\affil[1]{School of Engineering and Applied Sciences, Harvard University, Cambridge, MA 02134, USA}

\date{}
\maketitle            
\begin{center}
Submitted to Experimental Mathematics
\end{center}

\begin{abstract}

Based on the data of 12-17-crossing knots, we establish three new conjectures about the hyperbolic volume and knot cohomology:

(1) There exists a constant $a\in \mathds{R}_{>0}$ such that the percentage of knots for which the following inequality holds converges to 1 as the crossing number $c \to \infty$: 

 $$ \log r(K) < a \cdot \mbox{Vol}(K) $$ 
for a knot $K$ where $r(K)$ is the total rank of knot Floer homology (KFH) of $K$ and $\mbox{Vol}(K)$ is the hyperbolic volume of $K$.

(2) There exist constants $a,b\in \mathds{R}$ such that the percentage of knots for which the following inequality holds converges to 1 as the crossing number $c \to \infty$:

 $$ \log \mbox{det}(K) < a \cdot \mbox{Vol}(K) +b$$

  for a knot $K$ where $\mbox{det}(K) $ is the  knot determinant of $K$.

(3) Fix a small cut-off value $d$ of the total rank of KFH and let $f(x)$ be defined as the fraction of knots whose total rank of knot Floer homology is less than $d$ among the knots whose hyperbolic volume is less than $x$. Then for sufficiently large crossing numbers, the following inequality holds

 $$f(x) < \frac{L}{1 + \exp(-k\cdot(x-x_0))} + b$$

    where $L ,x_0, k, b$ are constants. 
\end{abstract}

\section{Introduction}
A major goal of knot theory is to relate the geometric structure of a knot complement to the knot’s topological properties. 
Nathan Dunfield \cite{dunfield} documented a nearly linear relationship between $\log \mbox{det(K)}$ and the hyperbolic volume of $K$ for all alternating knots $K$ with at most 13 crossings and samples of 14-16-crossing alternating knots ($\mbox{det(K)}$ denotes the knot determinant of $K$).  He then considered $\log \mbox{det(K)}/\log(\mbox{deg(J)})$ \textit{versus} the hyperbolic volume of $K$ for all alternating knots $K$ with at most 13 crossings and samples of 14-16-crossing alternating knots  ($\mbox{deg(J)}$ denotes the degree of the Jones polynomial of $K$).  Dunfield again observed the points to cluster around a straight line. The plots, however, become more scattered when including non-alternating knots. Furthermore, Mikhail Khovanov \cite{khovanov} observed a correlation between $\log(\mbox{rank }  H(K)- 1)$ and the hyperbolic volume of 10- and 11-crossing non-alternating knots where $\mbox{rank }H(K)$ is the total rank of Khovanov homology $H(K)$ of $K$.

In this work, we are interested in investigating the relationship between the hyperbolic and homological measures of the knot complexity. Recall that  the vast majority of knots are hyperbolic \cite{Thurston:1982zz} and that the volume of a knot complement serves as a proxy for the complexity of the knot. Another proxy of the complexity of a knot is the total rank of its knot Floer homology. In particular, the unknot is the only known knot whose total rank of knot Floer homology is equal to 1. In this paper, we quantify the strength of the nearly linear relation between the hyperbolic volume and the logarithm of the total rank of knot Floer  homology of knots with 12-17 crossings.
We then perform a similar analysis considering  the hyperbolic volume \textit{versus} the logarithm of the knot determinant. 
We find that the relationship between the hyperbolic volume and the  logarithm of the total rank of knot Floer homology is stronger compared to the relationship between the hyperbolic volume and the logarithm of the knot determinant. Finally, we provide experimental evidence for a special pattern in the density of knots with small total ranks of knot Floer homology.

This paper is organized as follows. Section \ref{section:background} provides a brief overview of the  concepts relevant to this work. In Section \ref{section:data}, we discuss how we constructed our dataset. In Section \ref{section:volume-and-kfh}, we provide experimental evidence for and formulate a conjecture that there exists an inequality between the hyperbolic volume and the logarithm of the total rank of knot Floer homology. In Section \ref{section:vol-and-determinant}, we conjecture the existence of an inequality between the hyperbolic volume and the logarithm of the knot determinant.  
Section \ref{section:density}  establishes a conjecture about the density of knots with small total ranks of knot Floer homology  based on our experimental data. Section \ref{section:conclusion} provides a summary of our findings. 

\section{Background information}\label{section:background}
Our primary focus in this paper is on the hyperbolic volume of the knot complement, knot Floer homology, and the knot determinant. We now briefly discuss each of these concepts. 
Knot Floer homology was introduced independently by Ozsváth-Szabó \cite{https://doi.org/10.48550/arxiv.math/0209056} and Rasmussen \cite{https://doi.org/10.48550/arxiv.math/0306378} around 2002. Knot Floer homology categorifies the Alexander polynomial and contains information about several non-trivial geometric properties of the knot (genus, slice genus, fiberedness, effects of surgery, and others). Knot Floer homology provides more information about knots/links compared to the Alexander polynomial. For example, it detects the genus of a knot while the Alexander polynomial gives only bounds for it. In addition, knot Floer homology detects fibered knots while the Alexander polynomial gives obstructions about it. Another advantage of knot Floer homology is that it is computable.  

Next, we define the hyperbolic volume of the complement of a hyperbolic link. A hyperbolic link is a link in the 3-sphere whose complement (the space formed by removing the link from the 3-sphere) can be given a complete Riemannian metric of constant negative curvature, giving it the structure of a hyperbolic 3-manifold, a quotient of hyperbolic space by a group acting freely and discontinuously on it. The components of the link will become cusps of the 3-manifold, and the manifold itself will have a finite volume. Mostow's Rigidity Theorem (see \cite{10.2307/j.ctt1k3s9kd}) implies that the hyperbolic structure on a finite volume hyperbolic 3-manifold is unique up to isometry. In particular, any invariants which are defined in terms of the hyperbolic structure of such a manifold are topological invariants of the 3-manifold. One of such invariants is the hyperbolic volume of the manifold. If the manifold has been decomposed into ideal hyperbolic tetrahedra, the volume will simply be the sum of the volumes of the tetrahedra \cite{Adams1991}.  The hyperbolic volume of the complement of a knot is a knot invariant.

The knot determinant has multiple equivalent definitions \cite{bams/1183539447}. The determinant $\mbox{det}(K)$ of a knot $K$ is defined as the absolute value of the determinant of the matrix $V +V^T$ where $V$ is a Seifert matrix of $K$. Equivalently, it can also be defined as $|\Delta_K(-1)|$ (i.e., the absolute value of the Alexander polynomial of $K$ evaluated at -1), which is the same as  $|J_{-1}(K)|$ (i.e., the absolute value of the Jones polynomial of $K$ evaluated at -1).  Finally, $\mbox{det}(K)$ is  equal to the number of elements in the first homology group of the double
cover of $S^3$ branched over $K$ \cite{Lickorish1997AnIT}.
\sloppy
\section{Data}\label{section:data} 
We used the knot database provided in \cite{burton:LIPIcs:2020:12183}, which contains all prime knots up to 19 crossings. In this work, we considered all knots with 12-17 crossings.  There are 1,288 alternating 	and 888 non-alternating  12-crossing knots; 4,877 alternating 	and 5,108 non-alternating 13-crossing knots; 19,536 alternating and 27,433 non-alternating 14-crossing knots; 85,262 alternating  and 168,023 non-alternating 15-crossing knots; 379,799 alternating  and 1,008,895 non-alternating 16-crossing knots; 1,769,978 alternating  and 6,283,385 non-alternating 17-crossing knots.

We converted the identification of the knots in this database to the Dowker–Thistlethwaite notation using the  \textsc{Regina} \cite{regina} Python package. We then used \textsc{Snappy} \cite{SnapPy} to compute the following invariants for each knot: the knot Floer homology,  the hyperbolic volume,  and the knot determinant.  
We have provided online\footnote{\url{doi.org/10.5281/zenodo.7879466}} the resulting dataset used in this work. 

\begin{figure*}[!htb]
\begin{center}
\includegraphics[width=1\textwidth]{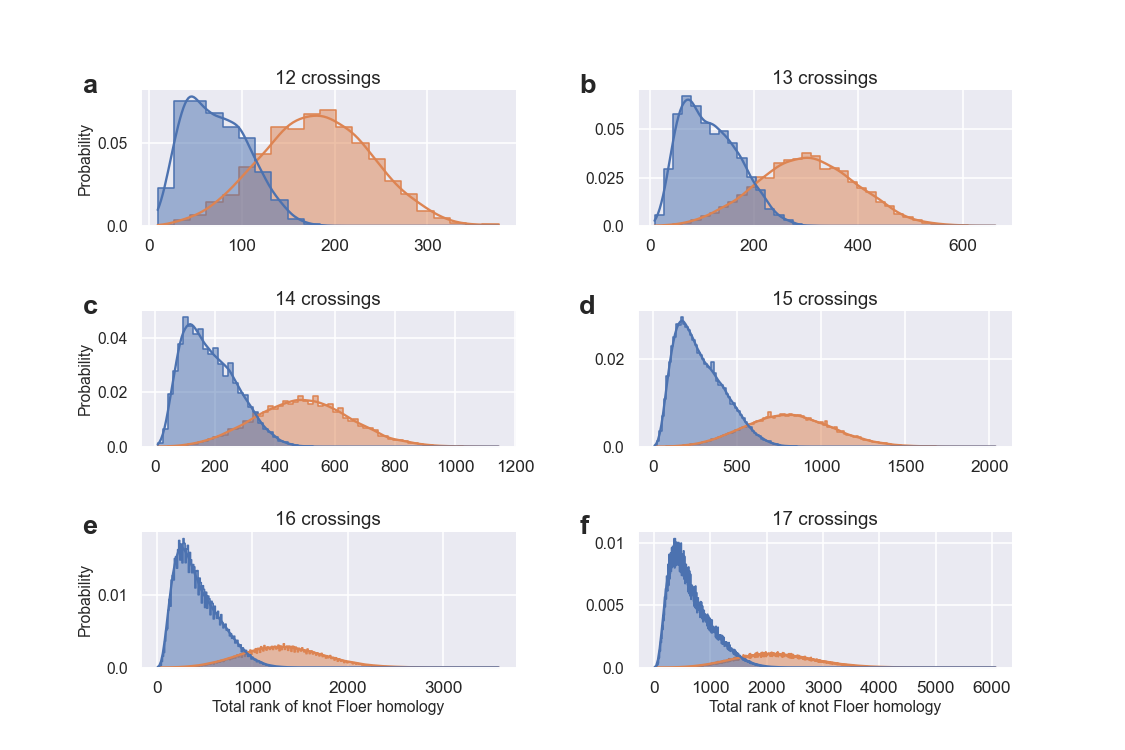}
\end{center}
\caption{The probability distribution of the total rank of knot Floer homology.  Non-alternating knots are shown in blue. Alternating knots are shown in orange.}
\label{fig:kfh-dist}
\end{figure*}

\begin{figure*}[!htb]
\begin{center}
\includegraphics[width=1\textwidth]{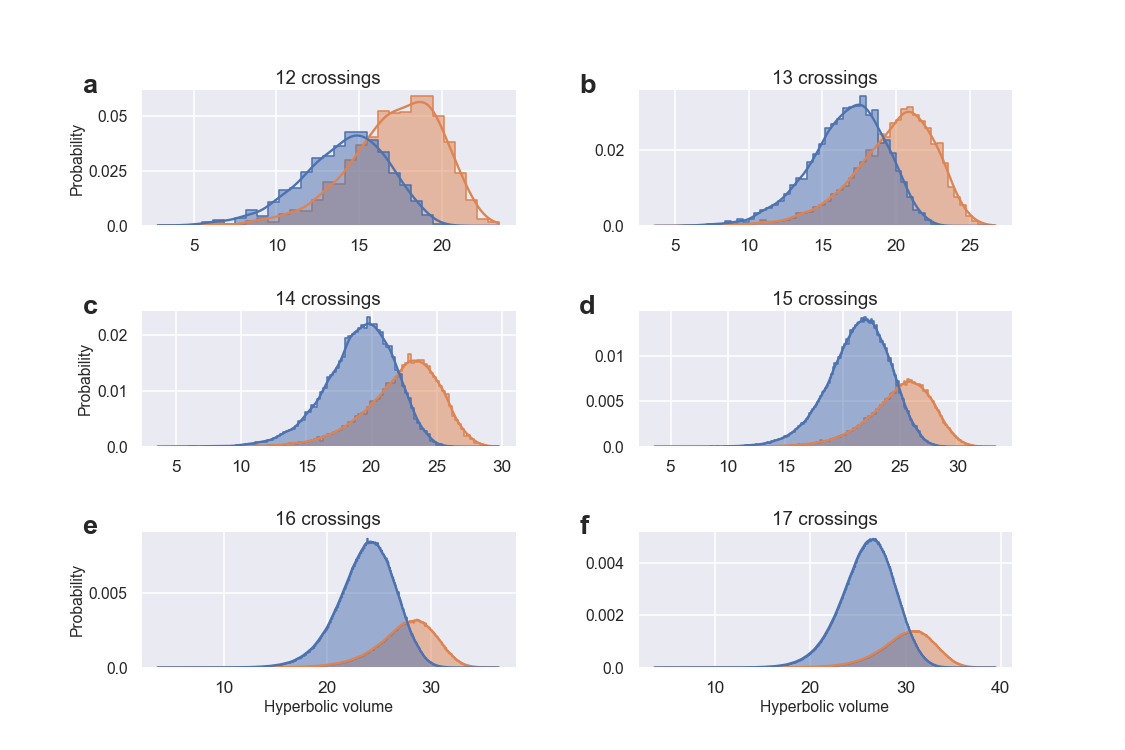}
\end{center}
\caption{The probability distribution of the hyperbolic volume.  Non-alternating knots are shown in blue. Alternating knots are shown in orange.}
\label{fig:volume-dist}
\end{figure*}

\begin{figure*}[!htb]
\begin{center}
\includegraphics[width=1\textwidth]{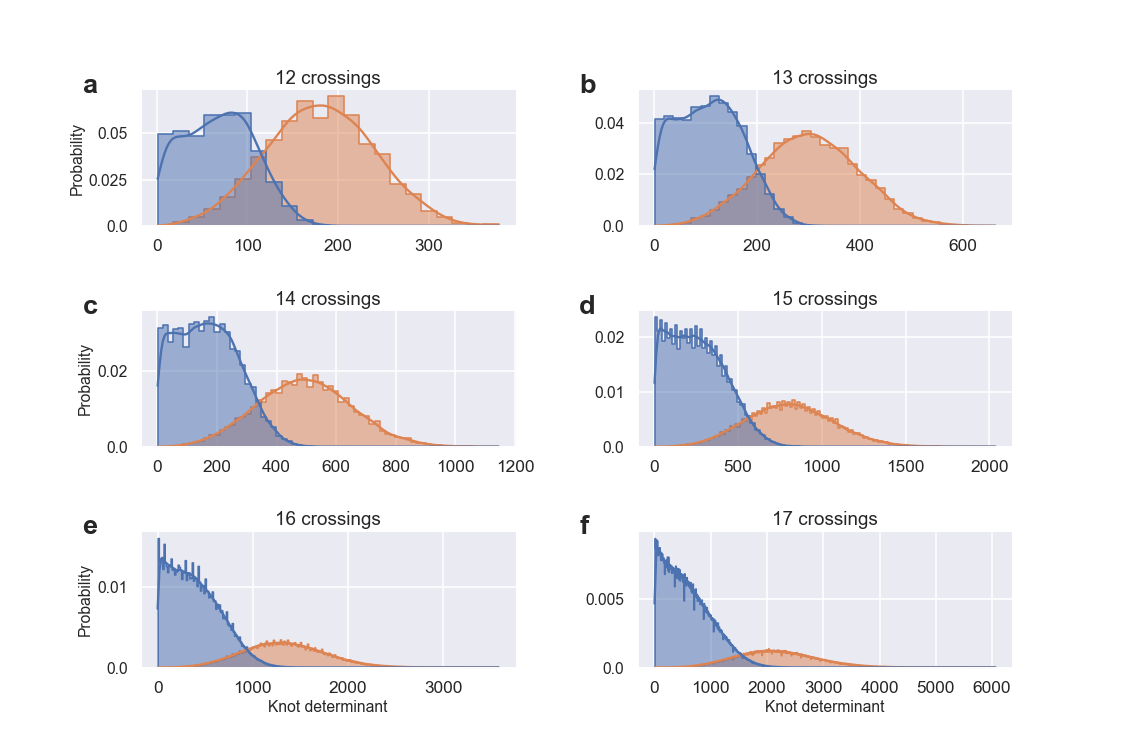}
\end{center}
\caption{The probability distribution of the knot determinant.  Non-alternating knots are shown in blue. Alternating knots are shown in orange.}
\label{fig:determinant-dist}
\end{figure*}

Figures \ref{fig:kfh-dist}, \ref{fig:volume-dist}, and \ref{fig:determinant-dist} show the probability distribution functions of the total rank of knot Floer homology, the hyperbolic volume,  and the  knot determinant, respectively, for 12-17-crossing knots.
Note that for a fixed crossing number, the mean value of a given knot invariant is always lower for the non-alternating knots compared to the alternating knots. 

\section{Volume and Cohomology}\label{section:volume-and-cohomology}
In all of our experiments, we considered knots of each crossing number (ranging between 12-17) and each type (alternating \textit{versus} non-alternating) separately.  We used least-squares linear regression to fit the hyperbolic volume \textit{versus} the following two invariants: the logarithm of the total rank of knot Floer homology and the logarithm of the knot determinant. The outputs of the linear regression fits included the estimated slope with its standard error, the intercept with its standard error, the Pearson product-moment correlation coefficient, and the $R^2$ value. We also plotted the hyperbolic volume \textit{versus} the number of diagonals that contain nontrivial knot Floer homology groups of a given knot. We did not observe any clear pattern in this case and did not include these results in the manuscript.

\subsection{Volume and Knot Floer Homology}\label{section:volume-and-kfh}

\begin{figure*}[!h]
\begin{center}
\includegraphics[width=1\textwidth]{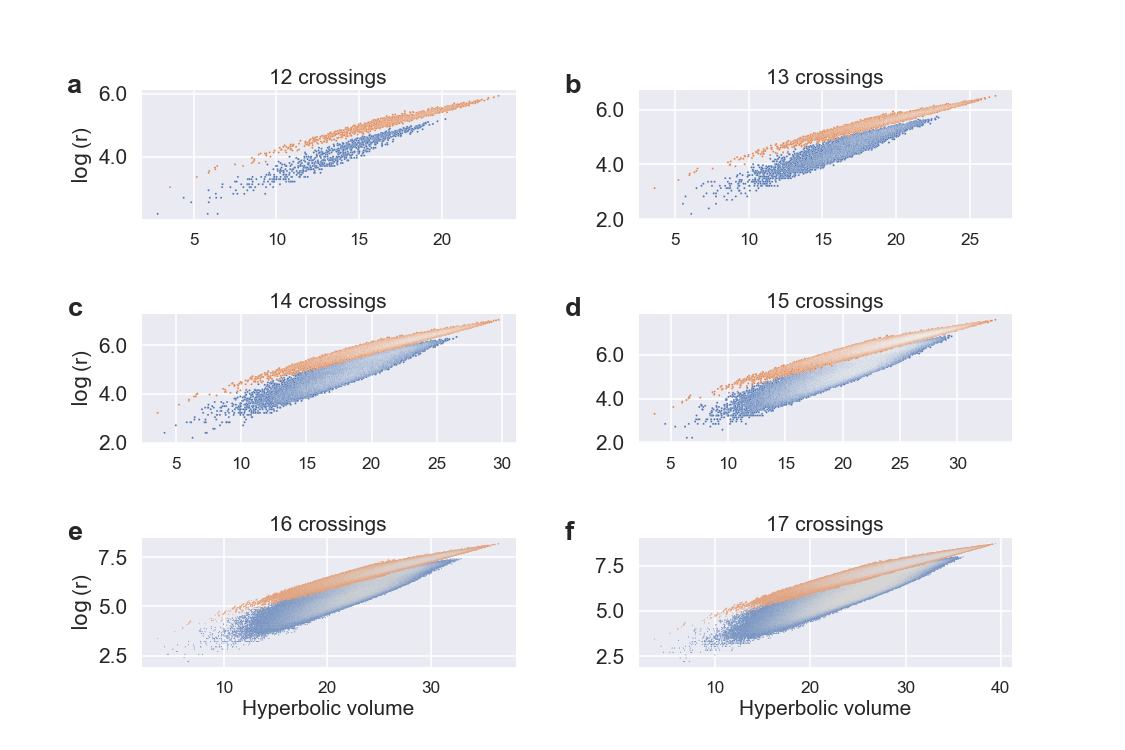}
\end{center}
\caption{Volume \textit{versus} the logarithm of the total rank of Knot Floer homology $r$.  Non-alternating knots are shown in blue. Alternating knots are shown in orange.}
\label{fig:log_kfh_rank_vs_volume}
\end{figure*}

Figure \ref{fig:log_kfh_rank_vs_volume} shows plots of the hyperbolic volume \textit{versus} the logarithm of the total rank of knot Floer homology  for 12-17-crossing knots. We see a linear trend in the case of alternating and non-alternating knots, although the data have more spread in the case of non-alternating knots. Table~\ref{table:fits-kfh-volume} summarizes the results of the linear regression fits of  the hyperbolic volume \textit{versus} the logarithm of the total rank of knot Floer homology for knots of a fixed crossing number (ranging between 12-17) and a fixed type (alternating \textit{versus} non-alternating). Overall, the correlation is strong in both cases, ranging between 0.97 and 0.982 for alternating knots and between 0.903 and 0.95  for non-alternating knots. The $R^2$ ranges between 0.941 and 0.964 for alternating knots and between 0.815 and 0.902 for non-alternating knots.

For alternating knots, the intercept value overall increases with increasing crossing number (except for 14-crossing knots), going from $2.94\pm 0.01$ for 12-crossing  knots to $3.989 \pm 0.0007$ for 17-crossing  knots. For non-alternating knots, on the other hand, the intercept value overall decreases with increasing crossing number (except for 14-crossing knots), going from $1.56\pm 0.03$ for 12-crossing knots to $0.124\pm  0.001$  for 17-crossing  knots. 

For alternating knots, the slope value decreases with increasing crossing number, going from $0.1284\pm 0.0007$ for 12-crossing  knots to $0.12083 \pm 2 \times 10^{-5}$ for 17-crossing  knots. For non-alternating knots, on the other hand, the slope value overall increases with increasing crossing number (except for 14-crossing knots), going from $0.185\pm 0.002$ for 12-crossing knots to $0.19507\pm  4\times 10^{-5}$  for 17-crossing  knots.

\begin{table}[!h]
\caption{Volume \textit{versus} the logarithm of the total rank of knot Floer homology.}
\begin{threeparttable}
\begin{tabular}{|l|*{7}{c|}}
\hline
$c$
&\makebox{$R^2$}&\makebox{Correlation}&\makebox{Slope}&\makebox{Slope error} &\makebox{Intercept}&\makebox{Intercept error}\\\hline

12n   & 0.902 & 0.95 & 0.185 & 0.002 & 1.56 & 0.03\\\hline

13n  & 0.872 & 0.934 & 0.189 & 0.001 & 1.49 & 0.02\\\hline

14n &  0.848 & 0.921& 0.1872 & 0.0005 & 1.507&0.009\\\hline

15n & 0.836 & 0.914&0.1913 & 0.0002 & 1.401 & 0.004\\\hline

16n & 0.821 & 0.906 & 0.19281 & 9 $\times 10^{-5}$&1.334&0.002\\\hline
17n &0.815&0.903&0.19507& $4\times 10^{-5}$&1.240&0.001\\\hline
12a  &0.964 & 0.982 & 0.1284 & 0.0007 & 2.94 & 0.01\\\hline
13a &0.958 & 0.979&0.1244&0.0004&3.183&0.007 \\\hline
14a &0.955 & 0.977&0.1239&0.0002&3.372&0.004\\\hline
15a &0.949& 0.974& 0.1221 &
$1\times 10^{-4}$&3.593&0.002\\\hline
16a &0.946&0.973& 0.12147&$5\times 10^{-5}$&3.791&0.001\\\hline
17a &0.941&0.97&0.12083&$2 \times 10^{-5}$&3.989&0.0007\\\hline
\end{tabular}
\begin{tablenotes}  
\item[*] $c$ denotes what data have been fitted, specified by the crossing number and  a label telling whether the knots are alternating (“a”) or non-alternating (“n”). $R^2$ is the coefficient of determination of the best-fit linear model. Correlation is the Pearson product-moment correlation coefficient between the hyperbolic volume and the total rank of knot Floer homology. Slope is the slope of the regression line.  Slope error is the standard error of the estimated slope under the assumption of residual normality. Intercept is the intercept of the regression line. Intercept error is the standard error of the estimated intercept under the assumption of residual normality.\label{table:fits-kfh-volume}
\end{tablenotes}
\end{threeparttable}
\end{table}

For each crossing number $c$ between 12-17, we also computed the minimum slope $a_{\mbox{min}}$ such that all of the knots in Figure \ref{fig:log_kfh_rank_vs_volume}  with crossing number $c$ lie below the line with the slope equal to $a_{\mbox{min}}$ and the intercept equal to 0. The results are shown in Table \ref{table:a_min}. $a_{\mbox{min}}$ slowly increases with increasing  crossing number, going from 0.852 for 12-crossing knots to 0.948 for 17-crossing knots. 

 \begin{table}[!h]
\begin{center}
\begin{tabular}{|l|*{7}{c|}}
\hline
$c$&\makebox{12}&\makebox{13}&\makebox{14}&\makebox{15}&\makebox{16}&\makebox{17}\\\hline

\makebox[3em]{$a_{\mbox{min}}$} & 0.852&0.874& 0.894&0.91309& 0.931&0.948\\\hline
 
\end{tabular}
\caption{$c$ is the crossing number. $a_{\mbox{min}}$ is the corresponding minimum slope  such that all of the knots in Figure \ref{fig:log_kfh_rank_vs_volume}  with  crossing number $c$ lie below the line with the slope equal to $a_{\mbox{min}}$ and the intercept equal to 0.\label{table:a_min}}
\end{center}
\end{table}

Our experimental data presented in this Section  led us to consider the following conjecture relating the total rank of knot Floer homology and the hyperbolic volume:

\begin{conjecture} \label{conjecture:volume-kfh} There exists a constant $a\in \mathds{R}_{>0}$ such that
  $$ \log r(K) < a \cdot Vol(K) $$

  for all knots $K$ where $r(K)$ is the total rank of knot Floer homology of $K$ and $Vol(K)$ is the hyperbolic volume of $K$.
\end{conjecture}

 However, conjecture \ref{conjecture:volume-kfh} does not hold in full generality. The following family of pretzel knots, noted by Nathan Dunfield, provides a counterexample. Let $K_n$ be the $(3, 3, n)$ pretzel knot, where $n$ is even.  Then $Vol(K_n)$ is uniformly bounded (less than 10, in fact), but $|\mbox{det}(K_n)|$ is $6n + 9$ which goes to infinity as $n$ goes to infinity. Since twist knots are alternating, the total rank of the knot Floer homology is equal to the determinant. Thus, Conjecture \ref{conjecture:volume-kfh} is not true.

Another counterexample is the following, which appeared in \cite{misev_spano_2021}  and was brought to our attention by 
David Futer. The authors of \cite{misev_spano_2021} build a sequence of fibered hyperbolic knots $K_n$, of fixed genus $g$, so that the rank of the knot Floer homology of $K_n$ in grading $(g-1)$ grows quadratically with $n$. At the same time, the monodromies of these knots are such that volume is universally bounded for each $g$. So the volume stays bounded but the total rank of the knot Floer homology of $K_n$ grows, contradicting the conjecture.

We therefore reformulate our conjecture as follows:

\begin{conjecture} \label{conjecture:volume-kfh-almost-surely} 
There exists a constant $a\in \mathds{R}_{>0}$ such that the percentage of knots for which the following inequality holds converges to 1 as the crossing number $c \to \infty$: 

 $$ \log r(K) < a \cdot \mbox{Vol}(K) $$ 
for a knot $K$ where $r(K)$ is the total rank of knot Floer homology (KFH) of $K$ and $\mbox{Vol}(K)$ is the hyperbolic volume of $K$. 
\end{conjecture}

In other words, the conjecture states that   $ \log r(K) < a \cdot Vol(K) $
holds asymptotically almost surely.

\subsection{Volume and Knot Determinant}  \label{section:vol-and-determinant} 

\begin{figure*}[!h]
\begin{center}
\includegraphics[width=1\textwidth]{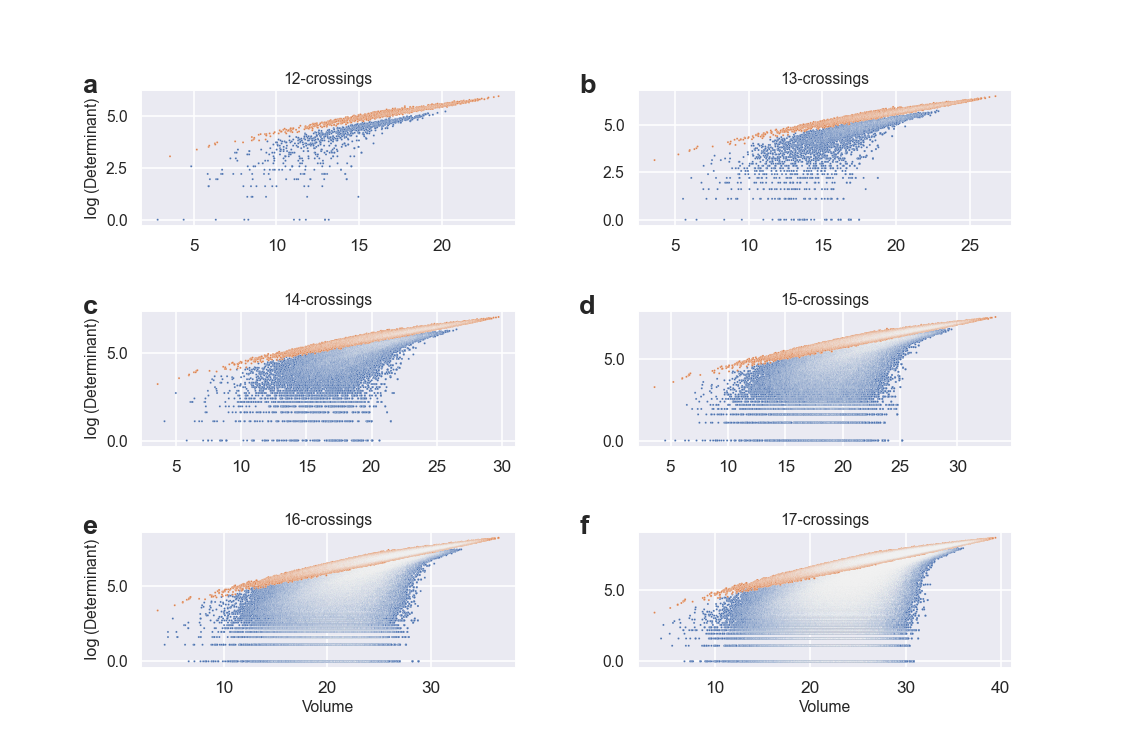}
\end{center}
\caption{Volume \textit{versus} the logarithm of the knot determinant.  Non-alternating knots are shown in blue. Alternating knots are shown in orange.}
\label{fig:log_det_vs_volume}
\end{figure*}

Figure \ref{fig:log_det_vs_volume} shows plots of  the hyperbolic volume  \textit{versus} the logarithm of the knot determinant for 12-17-crossing knots. Recall that the total rank of knot Floer homology is the same as the knot determinant for all alternating knots, so the data shown in Figures   \ref{fig:log_kfh_rank_vs_volume} and \ref{fig:log_det_vs_volume} coincide for alternating knots. 

Comparing the two Figures, we see that in  Figure \ref{fig:log_det_vs_volume}, the  data for the non-alternating knots have substantially more spread. Table~\ref{table:fits-det-volume} summarizes the results of the linear regression fits of the hyperbolic volume 
 \textit{versus} the logarithm of the knot determinant for non-alternating knots of a fixed crossing number ranging between 12 and 17. Overall, the $R^2$ (which ranges between 0.461 and 0.608) and the correlation (which ranges between 0.679 and 0.78) are much lower than the ones we observed for non-alternating knots in Section \ref{section:volume-and-kfh}, in which case $R^2$ ranged between 0.815 and 0.902 and the correlation ranged between 0.903 and 0.95. These metrics suggest that the hyperbolic volume of the knot complement correlates more strongly with the logarithm of the total rank of knot Floer homology compared to the logarithm of the knot determinant.

 \begin{table}[!h]
\caption{Volume \textit{versus} the logarithm of the knot determinant.}
\begin{threeparttable}
\begin{tabular}{|l|*{7}{c|}}
\hline
$c$
&\makebox{$R^2$}&\makebox{Correlation}&\makebox{Slope}&\makebox{Slope error} &\makebox{Intercept}&\makebox{Intercept error}\\\hline

12n  &0.608&0.78& 0.266 & 0.007 & 0.2 & 0.1\\\hline

13n  & 0.532& 0.729&0.271& 0.004&-0.13&0.06\\\hline

14n  & 0.499&0.706&0.263&0.002&-0.21&0.03\\\hline

15n & 0.484& 0.696 &0.2680& 0.0007&-0.53&0.01 \\\hline

16n & 0.465& 0.682& 0.2674&0.0003& -0.750&0.007\\\hline
17n &0.461 & 0.679 &0.2679&0.0001&-0.992&0.003\\\hline
\end{tabular}
\begin{tablenotes}
 \item[*] $c$ denotes what data have been fitted, specified by the crossing number and  a label telling whether the knots are alternating (“a”) or non-alternating (“n”). $R^2$ is the  coefficient of determination of the best-fit linear model. Correlation is the Pearson product-moment correlation coefficient between the hyperbolic volume and the total rank of knot Floer homology. Slope is the slope of the regression line.  Slope error is the standard error of the estimated slope under the assumption of residual normality. Intercept is the intercept of the regression line. Intercept error is the standard error of the estimated intercept under the assumption of residual normality.\label{table:fits-det-volume}
\end{tablenotes}
\end{threeparttable}
\end{table}

 Our experimental data presented in this Section  led us to consider the following conjecture relating the knot determinant and the hyperbolic volume:

\begin{conjecture} \label{conjecture:volume-det} There exist constants $a,b\in \mathds{R}$ such that
  $$\log \mbox{det}(K) < a \cdot Vol(K) +b$$

  for all knots $K$ where $\mbox{det}(K) $ is the determinant of $K$ and $Vol(K)$ is the hyperbolic volume of $K$.
\end{conjecture}


However, the conjecture does not hold in full generality. The following family of pretzel knots, noted by Nathan Dunfield, provides a counterexample. Let $K_n$ be the $(3, 3, n)$ pretzel knot, where $n$ is even.  Then $Vol(K_n)$ is uniformly bounded (less than 10, in fact), but $|\mbox{det}(K_n)|$ is $6n + 9$ which goes to infinity as $n$ goes to infinity. 

In general, as Andras Juhasz noted,  for any knot, if you start twisting two strands, the volume converges (this is a special case of Thurston's hyperbolic Dehn surgery theorem). In particular, the volume of any twist knot is bounded. On the other hand, the determinant is $2n+1$ for an $n$-twist knot. See \cite{Stoimenow2007GRAPHSDO} for a more detailed discussion of the relation between the volume and the knot determinant. \cite{Stoimenow2007GRAPHSDO} proved that if $L$ is a nontrivial nonsplit alternating link, then $\mbox{det}(L) \geq 2 \cdot 1.0355^{\mbox{Vol(L)}}$. \cite{Stoimenow2007GRAPHSDO}  also proved that there are constants $C_1, C_2 > 0$ such that, for any hyperbolic alternating link $L$, $\mbox{det}(L) \leq (\frac{C_1 \cdot c(L)}{\mbox{Vol(L)}})^{C_2 \mbox{Vol(L)}}$ where $c(L)$ is the crossing number of a link $L$.

We therefore reformulate our conjecture as follows:

\begin{conjecture} \label{conjecture:volume-det-almost-surely} 
There exist constants $a,b\in \mathds{R}$ such that the percentage of knots for which the following inequality holds converges to 1 as the crossing number $c \to \infty$:

 $$ \log \mbox{det}(K) < a \cdot \mbox{Vol}(K) +b$$

  for a knot $K$ where $\mbox{det}(K) $ is the  knot determinant of $K$.

\end{conjecture}

 \subsection{Knot Floer Homology Total Rank Density}\label{section:density}

\begin{figure*}[!h]
\begin{center}
\includegraphics[width=1\textwidth]{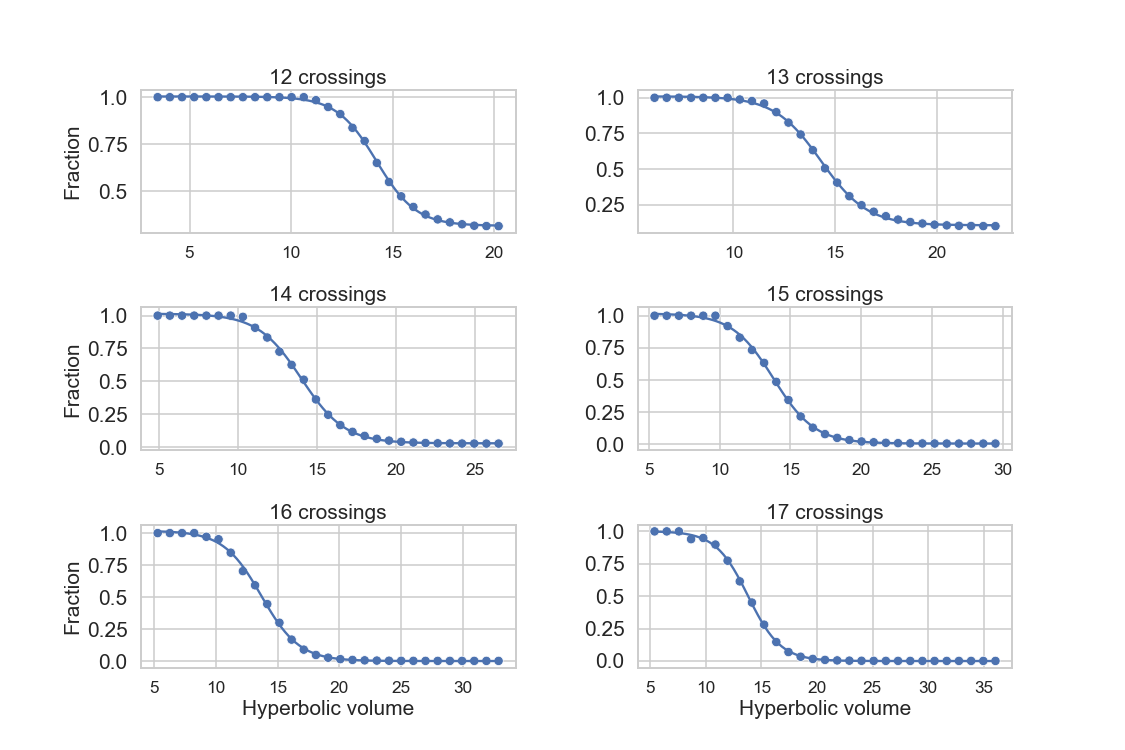}
\end{center}
\caption{Fraction of non-alternating knots whose total rank of knot Floer homology is less than $d=50$ among the non-alternating knots whose hyperbolic volume is less than a given value. The line shows the best-fit sigmoid function.}
\label{fig:dist_50}
\end{figure*}

Non-alternating knots with small total ranks  of knot Floer homology usually have very special properties (e.g., they can be torus knots) and are, therefore, of mathematical interest.  We performed the following analysis to study the distribution of non-alternating knots with small total ranks of knot Floer homology.  We fixed the cut-off value of the total rank of knot Floer homology to be $d$ (say, 50). Next, we ordered the non-alternating knots of a fixed crossing number by the hyperbolic volume and considered all knots whose hyperbolic volume is less than a given value $x$. Among these knots, we then computed the fraction of knots whose  total rank of knot Floer homology was less than $d$. We called the function defined in this way $f(x)$. Figure \ref{fig:dist_50} shows a plot of this function when $d=50$. 

Having observed a consistent sigmoid-like shape of $f(x)$  for all non-alternating knots with the  crossing number ranging between 12-17, we fitted the data with a sigmoid function defined as 
\begin{equation}\label{eq:sigmoid}
g(x) = \frac{L}{1 + \exp(-k\cdot(x-x_0))} + b.
\end{equation}

The best-fit parameters are reported in Table \ref{table:density1}. The results of our fits lead us to conjecture that such density can be bounded by a sigmoid function in the limit of large crossing numbers. We summarize this finding in the following conjecture:  

\begin{table}[!h]
\caption{Knot Floer Homology Total Rank Density Function. $d=50$.}
\begin{threeparttable}
\begin{tabular}{|l|*{9}{c|}}
\hline
$c$& \makebox{$L$}&\makebox{$L$ error}&\makebox{$k$}&\makebox{$k$ error}&\makebox{$x_0$}&\makebox{$x_0$ error}&\makebox{$b$}&\makebox{$b$ error}\\\hline

12n&-0.689& 0.003 & 14.21	& 0.02 &1.02&0.02 & 1.004& 0.002\\\hline

13n&-0.902&0.005	&14.31&0.02 &0.88&	0.02&1.01 & 0.003\\\hline

14n&-0.986 &0.007 &	14.00	&0.04 &0.73&	0.02 & 1.014& 0.005\\\hline

15n&-1.01& 0.007 &	13.80	& 0.04 &0.69&0.02 &	1.017 &0.006 \\\hline

16n &-1.017&0.007	&13.69& 0.05 &0.65& 0.02	& 1.017&  0.006\\\hline

17n &-1.002& 0.006 &	13.81	& 0.04&0.68& 0.01&	 1.002&0.005 \\\hline
\end{tabular}
\begin{tablenotes}
\item[*] $c$ denotes what data have been fitted, specified by the crossing number and a label telling that the knots are non-alternating (“n”).  $L,k,x_0$, and $b$ are the parameters of the best-fit sigmoid model defined in equation \ref{eq:sigmoid}. $L$ error, $k$ error, $x_0$ error, and $b$ error are the one standard deviation errors of $L,k,x_0$, and $b$ parameters, respectively. \label{table:density1}
\end{tablenotes}
\end{threeparttable}
\end{table}

\begin{conjecture} \label{conjecture:density}
Fix a  small cut-off value $d$ of the  total rank of knot Floer homology. Let $f(x)$ be defined as the fraction of knots whose total rank of knot Floer homology is less than $d$ among the knots whose hyperbolic volume is less than $x$. Then for sufficiently large crossing numbers, the following inequality must hold

$$f(x) < \frac{L}{1 + \exp(-k\cdot(x-x_0))} + b$$

    where $L ,x_0, k, b$ are constants.
\end{conjecture}

\section{Conclusion}\label{section:conclusion}

 By analyzing data of 12-17-crossing knots, we experimentally examined and quantified the strength of the following two knot invariant relationships: (i) hyperbolic volume \textit{versus} the logarithm of the total rank of knot Floer homology and (ii) hyperbolic volume \textit{versus} the logarithm of the knot determinant. We concluded that relationship (i) is stronger than relationship (ii), with substantially higher $R^2$ and correlation values. Based on our experimental findings, we formulated a new conjecture that there exists an inequality between the hyperbolic volume and the  logarithm of the total rank of knot Floer homology which holds asymptotically almost surely (Conjecture \ref{conjecture:volume-kfh-almost-surely}). We also conjecture that there is an  inequality between the hyperbolic volume and the logarithm of the knot determinant which holds asymptotically almost surely (Conjecture \ref{conjecture:volume-det-almost-surely}). Further, having carried out computational experiments, we  discovered Conjecture  \ref{conjecture:density} that, for sufficiently large crossing numbers,  provides a sigmoid bound for the density of  non-alternating knots with small total ranks of knot Floer homology.  

We have  provided online the dataset\footnote{\url{doi.org/10.5281/zenodo.7879466}} and code\footnote{\url{https://github.com/eivshina/patterns-in-knot-floer-homology}} used in this work. We hope that these resources will be useful to anyone interested in further studying and understanding knot Floer homology and the relationships between different knot invariants.

\section{Code and dataset}
The code is available on GitHub at \url{github.com/eivshina/patterns-in-knot-floer-homology}. The dataset is archived on Zenodo at \url{doi.org/10.5281/zenodo.7879466}.

\section{Acknowledgments}
Thank you to Zoltán Szabó for advising this work. Thank you to Peter Ozsváth, Ian Zemke, Nathan Dunfield, Andras Juhasz, David Futer, and Tetsuya Ito for providing helpful feedback on the manuscript.

\bibliography{main}

\end{document}